\documentclass[12pt,a4paper]{article}
{\fontsize{12}{24}\selectfont 
\RequirePackage{amsmath}
\RequirePackage{relsize}
\RequirePackage{amsthm} 
\RequirePackage{amsfonts}
\RequirePackage{amssymb}
\RequirePackage{url}
\RequirePackage{physics}

\RequirePackage[bottom=3cm, right=2.5cm, left=2.5cm, top=3cm]{geometry}

\ProvideTextCommand{\DJ}{OT1}{\leavevmode\raisebox{-.5ex}{\makebox[0pt][l]{\hskip-.07em\accent"16\hss}}D}

\newtheorem{remark}{Remark}
\newtheorem{lemma}{Lemma}
\newtheorem{corollary}{Corollary}
\newtheorem{theorem}{Theorem}
\newtheorem{example}{Example}

\author{Alexander Kushkuley \\ kushkuley@gmail.com}
\title{ Some Remarks on Commuting Probability}
\begin{document}
\maketitle

\begin{abstract}
\noindent  
We introduce a weighted sum of irreducible character ratios as an estimator for commutator probabilities. The estimator yields Frobenius formula when  applied to a regular representation

\end{abstract}

\section{Introduction}
  \numberwithin{equation}{section}
  Let $G$ be  a finite group. 
The probability of $ g \in G $ to be a commutator 
(cf. e.g. \cite{Howe}, \cite{Serre2}) is defined as
\begin{align}
	c(g)  = \frac{| \{ [x,y] = g, \; x,y \in G  \}| } { | G\times G |} \nonumber
\end{align}  
\noindent The commuting probability (cf. \cite{Howe}-\cite{Gustafson}) of $ G $ can be defined thereby  as $ c(1) $
that is
\begin{align}
 	c(G) = c(1)  = \frac{| \{ [x,y] = 1, \; x,y \in G  \} |} { | G\times G |}  \nonumber
\end{align}

\noindent All groups considered in this paper are finite and all representations are assumed to be complex unitary. 
Let  $ \textnormal{Irr}(G) = \{ \psi_1, \cdots, \psi_k \} $ be  the set of all pairwise inequivalent complex characters of $ G, \; \dim \psi_i = \psi_i(1) = n_i, \; i = 1, \cdots , k $.

\bigskip \noindent First of all,  recall a well known formula of Frobenius (cf. \cite{Howe}, \cite{Serre2}) 
\begin{align}
	c(  g )	
	\; = \;  \frac{ (1/n_1) \psi_1(g) + \cdots + (1/n_k)  \psi_k(g) } 
	{|G|} 
\end{align}
and its immediate consequence (cf. \cite{Erdos})
\begin{align}
	c(G) = c(1) = k/|G| 
\end{align} 
\noindent Estimates of $ c(g)$ that we will introduce below are similar to (1.1), but in general, they do not involve all characters of $ G $. As we will see in Section 4, just one exact irreducible character is sufficient for a rough estimate of $c(g)$.  
\newline\newline\noindent Let $  m =|G| $ be the order of  $G$. It is assumed that $ G \times G $ is a probability space with standard (point mass $ 1/m^2$) Haar measure.   The expectation of  a random variable $x$ will be denoted by  $\mathbb{E}(x)$.
 
 \section{Statement of the result}
 We will say that a virtual character 
 \begin{align}
   \phi = \sum_{\chi \in \textnormal{Irr(G)}}{a_{\chi}} \chi, \;\; a_{\chi}\in \mathbb{C}   
 \end{align}
 of the group $G$  is non-negative if the following conditions hold:
  \begin{itemize}
  	\item[(a)] $ \Re(\phi(g) ) \geq 0 $ for all $ g \in G $ 
  	\item[(b)] all coefficients $ a_{\chi}  $ are non-negative real numbers and at least one of them is non-zero   
  \end{itemize}

\noindent  The kernel of a virtual character can be defined as
\begin{align}
	\ker (\phi ) = \{ g \in G, \; \phi(g) = \phi(1) \; \Longleftrightarrow \; \sum{a_i} \chi_i(g) = \sum{a_i} \chi_i(1) \; \} \nonumber
\end{align}

\noindent and we will say that a non-negative virtual character $ \phi $ is exact (faithful) if
$\ker(\phi) = \{ 1 \}$.  	 
   
 \begin{remark}
 Any permutation character is non-negative. Another example of a non-negative character is  $ n \cdot 1 + \chi $ for any character $ \chi $ of dimension $n$.  Any non-trivial combination of non-negative characters with non-negative coefficients is non-negative and it is exact if one of the summands is. The character $ \psi_1 + \cdots + \psi_k $ (Gelfand character of $G$) is known to be  non-negative for complex representations of ordinary finite classical groups, although it is not non-negative in general (see a discussion in \cite{Jorge}) 
 \end{remark}

 \noindent 
 For a given  virtual character 
 \begin{align}
  \phi = a_1 \chi_1 + \cdots + \ a_q \chi_q \; \; ( a_i \in \mathbb{R}, \; \chi_i \in \textnormal{Irr}(G), \; i = 1, \cdots, q) 
 \end{align}
 and a fixed $ g \in G $ define a real-valued random variable $\xi_{g,\phi}$ on $ G \times G $ as follows
 \begin{align}
 	\xi_{g,\phi} \equiv \xi_{g,\phi}(a,b) = \Re(\phi([a,b]^{-1}g) ) =  \sum_{i=1}^{q} a_i\Re(\chi_i([a,b]^{-1}g)), \; a,b \in G 
 \end{align} 	
\begin{remark}
	If $\phi $ is real valued then 
	 \begin{align}
	 \xi_{g,\phi}(a,b) = \phi([a,b]^{-1}g)  =  \sum_{i=1}^{q} a_i \chi_i([a,b]^{-1}g), \; a,b \in G  \nonumber
	\end{align} 	
\end{remark}
\noindent The following obvious lemma  illustrates these definitions
\begin{lemma} If the virtual character $\phi$ (2.2) is non-negative then  
\begin{align}	
c(g) \; \leq \; \textnormal{Prob}(\; [a,b]^{-1} g \in \ker ( \phi)  \; ) \; = \;  \textnormal{Prob}(\; \xi_{g,\phi} \; = \; \phi(1) \;   )  \nonumber
\end{align}
and if $ \phi $ is exact then 
\begin{align}	
	c(g) \; = \;  \textnormal{Prob}(\; \xi_{g,\phi} = \phi(1) \;   )  \nonumber
\end{align}
\end{lemma}

\noindent
We  can now formulate the main result of this paper. 

 \begin{theorem}
  	If $ \phi = a_1 \chi_1 + \cdots + \ a_q \chi_q \; (a_i \in \mathbb{R}, \; \chi_i \in \textnormal{Irr}(G), \; i = 1, \cdots, q) $
 	is a  virtual character of  $ G $  then
 \begin{align}
		\mathbb{E}(\xi_{g,\phi}) \; = \;   \sum_{i=1}^{q} (a_i/n_i^2) \Re ( \chi_i (g))
	\end{align}
If $\phi$ is non-negative then   for any  $ g \in G $   
\begin{align}
	c(  g )	\; \leq \; 	\mathbb{E}(\xi_{g,\phi}) / \phi(1) \; = \; 
 \frac{ (a_1/n_1^2) \Re(\chi_1(g)) + \cdots + (a_q/n_q^{2})  \Re(\chi_q(g)) } 
	{a_1  n_1 + \cdots + a_q  n_q}
\end{align}
\noindent and in particular
\begin{align}
	c(G) = c(1) \; \leq \;  \frac{(a_1/n_1) + \cdots + (a_q/n_q) }{a_1 n_1 + \cdots + a_q n_q }
\end{align}
\noindent where $ n_i =  \chi_i(1), \; i = 1, \cdots, q$
 \end{theorem}
 \noindent We precede the proof  by some auxiliary remarks 
 \section{Counting Lemmas (cf. \cite{Serre2})}  
 Let 
 $ \rho : G \rightarrow \textnormal{GL}(V) $ be an $n$-dimensional   representation of $ G $. Denote the character of $ \rho $ by $ \chi_{\rho}  $. For any  $ a \in \textnormal{End}(V) $ set 
 \begin{align}
 	A_{\rho} (a) = (1/m)\sum_{g \in G } \rho(g) a \rho(g)^{-1} \nonumber
 \end{align}
 Let also $ I = I_V $ denote the unity matrix in $\textnormal{End}(V) $   	
\begin{lemma} If representation $ \rho $ is irreducible then
for any $ a \in \textnormal{End}(V) $
	\begin{align}
		A_{\rho}(a) = (1/n) \textnormal{tr} (a) I_V  \nonumber
	\end{align} 
In particular,	
	if $h \in G $ then 
	\begin{align}
		A_{\rho}(h) = (1/n) \chi_{\rho}(h) I_V  \nonumber
	\end{align} 
\end{lemma}
 \noindent \textbf{Proof.}  $ 	A_{\rho} (a) $  commutes with $\rho(G)$ and therefore  $ 	A_{\rho} (a) = \lambda I_V $ for some $ \lambda \in \mathbb{C}$. 
 To find $\lambda$, note that
 \begin{align}
 	\textnormal{tr}(\lambda I_V ) = 	n \lambda = \textnormal{tr}(\; 		A_{\rho} (a)  \; ) =   \textnormal{tr}( a )   \nonumber  
 \end{align}
 and therefore $ \lambda = (1/n)\textnormal{tr}(a) $ as stated.
\newline\newline\noindent
For any representation $ \rho$ of $ G $ set  
\begin{align}
C_{\rho} = \frac{1}{m^2}\sum_{x,y\in G } \rho([x,y])  \nonumber
\end{align}  
and set 
\begin{align}
	 T_{\phi} = \frac{1}{m^2}\sum_{x,y\in G } \phi([x,y]) \nonumber
\end{align}
\noindent for any virtual character $\phi \in R_{\mathbb{C}}(G) $   

\begin{lemma} (cf. \cite{Serre2}, 8.13, Ex. 27). 
If representation $ \rho$ is irreducible then 
 $ C_{\rho} = (1/n)^2I$ and, therefore,
	$T_{\rho} = 1/n$ 
 
\end{lemma}
\noindent \textbf{Proof.} By Lemma 1,  

$$C_{\rho} = \frac{1}{m}\sum_{x \in G } A_{\rho}(x^{-1}) \rho(x) = 
\frac{1}{n} \left( \frac{1}{m}\sum_{x \in G } \chi_{\rho}(x^{-1}) \rho(x) \right) = (1/n^2)I
$$ 
As an obvious corollary, we get also
\begin{lemma} Let $\phi = a_1 \chi_1 + \cdots + a_q  \chi_q $ where $ \chi_i $ are irreducible characters of $ G $ and $ a_i \in \mathbb{C} $ are complex numbers, $ i = 1, \cdots , q $.  Then $$ T_\phi = \sum_{i=1}^q \frac{a_i}{n_i} $$
\end{lemma}

\section{Proof of Theorem 1, Corollaries and Examples}
\noindent  Using linearity and Lemma 2, we can establish formula (2.4) by direct computation
\begin{align}
	\!\!\!\!\ \mathbb{E}(\xi_{g,\phi}) \; = \; \sum_{i=1}^{q} a_i \frac{1}{m^2}\sum_{a,b \in G}\Re(\chi_i([a,b]^{-1}g)) \; = \;   
	\sum_{i=1}^{q} a_i \Re \left(\chi_i \left(\frac{1}{m^2}\sum_{a,b \in G}[a,b]^{-1} \cdot g \right) \right) \; = \; \nonumber \\
	= \;   \sum_{i=1}^{q} a_i \Re \left(\chi_i \left((1/n_i^2) \rho_i ( g ) \right) \right) \; = \nonumber  \; \sum_{i=1}^{q} (a_i/n_i^2) \Re ( \chi_i (g))
\end{align}
\noindent If $ \phi$  is non-negative then   
 $ \xi_{g,\phi} $ and its expectation  $ \mathbb{E}(\xi_{g,\phi}) $ are also non-negative.  In this case it is easy to check that  
$ \xi_{g,\phi}$ (2.3)  attains its maximum value $ \phi(1) = a_1 n_1 + \cdots + a_q n_q $ when (and only when)  $g = [a,b] $ for some $ a,b \in G $. Therefore, the estimate (2.5) follows from Markov inequality and established formula for the expectation of $\xi_{g,\phi}$  (2.4)    
\begin{align}
	c(g) \; \leq \; \textnormal{Prob}( \; \xi_{g,\phi} \; \geq \phi(1) \;  \; ) \; \leq \; \frac{\mathbb{E}(\xi_{g,\phi})} {\phi(1)}    \nonumber 
\end{align}	 
\noindent   The estimate (2.6)  is  a specialization of (2.5) for $ g = 1$. Note, however, that inequality (2.6) can be  established by direct computation of the expectation of a  random variable     
 $$ \xi \equiv \xi(a,b) = \sum_{i=1}^{q} a_i\Re(\chi_i([a,b])), \; a,b \in G $$
 (cf. Lemma 4)
 \newline\newline\noindent
 To avoid confusion we state the obvious
 \begin{corollary} (cf. Remark 2).
 	If in conditions of Theorem 1, the virtual character $\phi$ is real valued then
 	 \begin{align}
 		\mathbb{E}(\xi_{g,\phi}) \; = \;   \sum_{i=1}^{q} (a_i/n_i^2)  \chi_i (g)  \nonumber
 	\end{align}
 	and  if $ \phi $ is real valued and non-negative then for any  $ g \in G $   
 	\begin{align}
 		c(  g )	\; \leq \; 	\mathbb{E}(\xi_{g,\phi}) / \phi(1) \; = \; 
 		\frac{ (a_1/n_1^2) \chi_1(g) + \cdots + (a_q/n_q^{2})  \chi_q(g) } 
 		{a_1  n_1 + \cdots + a_q  n_q} \nonumber
 	\end{align}
 \end{corollary}  
\bigskip\noindent Let  
\begin{align}
	\mathfrak{r} = n_1 \psi_1 + \cdots + n_k \psi_k
\end{align} 
\noindent be the (real valued, exact and non-negative) regular character of $ G $ (cf.  Introduction).  Corollary 1 applied to $\mathfrak{r}$ reads   
 \begin{align}
 c(  g )	
	\; \leq \;  \frac{ (1/n_1) \psi_1(g) + \cdots + (1/n_k)\psi_k(g) }  	 
	{ |G|}  \;  = \;  \mathbb{E}_{g,\mathfrak{r}} / |G|, \;\; g \in G  \nonumber
\end{align}
\noindent It is obvious that   $ \mathbb{E}_{g,\mathfrak{r}} = |G| \cdot c(g) $ (cf. \cite{Serre2} or Lemma 1) and hence, formally speaking,  formula (1.1)  follows from Corollary 1. 
\newline\newline\noindent The vector $ ( n_1 = \psi_1(1), \cdots , n_k = \psi_k(1)  )$ of dimensions of irreducible characters of $G$ is a barycenter of an affine simplex in $  \mathbb{R}^k $ defined by  conditions

\begin{align}
	\sum_{i=1}^{k} \alpha_i n_i  =  |G|, \; \; 	\alpha_i \geq 0, \; i = 1, \cdots , k   	
\end{align}
\noindent and next corollary shows that  regular representation (4.1)  is an equilibrium point  in the space of non-negative virtual characters. Define the  $ L_1 $ norm of a complex function $f$ on the group $ G $ as  $ |f(g) |_1 = \sum_{g \in G } |\Re(\phi(g))| $
and consider the following  optimization problem
\begin{align}	
\underset{\alpha_1, \;\cdots, \;\alpha_k} {\textnormal{minimize}} \; \; |\sum_{i=1}^k (\alpha_i/n_i^{2}) \Re(\psi_i)|_1
\end{align}
\noindent under constraints (4.2) and
\begin{align}
 \sum_{i=1}^{k} \alpha_i \Re(\psi_i(g))  \geq 0, \; \; \textnormal{ for all } g \in G 
\end{align}
\begin{corollary} The regular character (4.1)
 	 is a  solution of the constrained minimization problem (4.2)-(4.4). 
\end{corollary}
\noindent \textbf{Proof}. Due to  constraints (4.4) and (4.2), the minimization takes place over non-negative virtual characters $ \phi$ . By Theorem 1 (2.4), the function that is minimized is non-negative. Hence, the norm sign can be removed from (4.1) and that allows to apply Theorem 1 (2.5) and the constraint (4.2) to finish the proof.

 \begin{corollary}
 	If $G$ has an irreducible representation of dimension $ n $ then 
 \begin{align}  c(g) \; \leq \; 
 	\frac{1}{2} \left( 1 + \frac{1}{n^3}\textnormal{Re}(\chi(g)) \right)  \end{align}
\noindent and in particular 	
 	\begin{align}
 		c(G) \; \leq \; \frac{1}{2} \left( 1 + \frac{1}{n^2} \right)   \tag{4.5'}
 	\end{align}
 \end{corollary}
 \noindent \textbf{Proof.} Let $ \chi \in \textnormal{Irr}(G) $ be an irreducible character of $ G $ of dimension $n$ ($\chi(1) = n$). Applying Theorem 1 to the non-negative  character $ n \cdot 1 + \chi$ (cf. Remark 1), we get (4.5). 

 \begin{example} (See \cite{Gustafson}).
  If $c(G) > 5/8$ then by (4.5') any irrep of $ G$ is one-dimensional and $G$ must be abelian in complete agreement with the "5/8 theorem" of Gustafson  (\cite{Gustafson}, \cite{comprob}). For the same reason, if $ c(G) = 5/8$ then all irreps of $ G $ must have dimension no greater than two.  It is well known  (cf. e.g. \cite{Gustafson}, \cite{comprob})
 	that commuting probability of $5/8$ is maximal for non-abelian  groups and this maximum  is attained, for example, by the group of quaternions that is a subgroup of   $ \textnormal{U}(2)$

 \end{example}
 
 \begin{corollary} (cf. Remark 1). The (permutational) character of the standard $n$-dimensional representation of the symmetric group $ S_n$ is equal to $ 1 + \chi_{n-1} $ where $  \chi_{n-1} $ is an irreducible chracter of dimension $ n-1$. Hence, we have estimates  
 	$$ c(g) \;  \leq \;  \frac{1}{n}\left(1 + \frac{ \chi_{n-1}(g) }{(n-1)^2} \right) \!, \: g \in S_n 
  $$ 
  and $$ c(S_n) = 1/(n-1) $$

 \end{corollary}
 
 \begin{corollary}
 	Let $\chi$ be an irreducible  character of  $G$ and let $$ \chi \otimes \bar{\chi} =  a_1 \chi_1 + \cdots + \ a_q \chi_q $$ where "Clebsch–Gordan coefficients" $a_i$ are positive integers  and  $ \chi_i$ are irreducible characters of $G$ of dimensions $ n_i,  \; i = 1, \cdots q  $. Then for any $ g \in G $
 	\begin{align}
 		c(  g )	
 		\; \leq \;  \frac{ (a_1/n_1^2) \chi_1(g) + \cdots + (a_q/n_q^{2})  \chi_q(g) } 
 		{a_1  n_1 + \cdots + a_q  n_q}  
 	\end{align}
 	and therefore 
 	\begin{align}
 		c(G) \; \leq \;  \frac{a_1/n_1 + \cdots + a_q/n_q }{a_1 n_1 + \cdots + a_q n_q } 
 	\end{align}
 \end{corollary}
 
 \begin{corollary}
 	Clebsch–Gordan coefficients $a_1, \cdots, a_q $ of an irreducible complex representation that is realizable over reals  must satisfy conditions (4.6), (4.7)  
 \end{corollary}
 
 \begin{example}
 	Let $\chi_5$ denotes the $5$-dimensional irreducible character of the alternating group $ A_5 $. All representations of $A_5$    are realizable over reals and it is easy to verify using the character table (cf. \cite{Serre2}) that
 	$$ \chi_5 \otimes \chi_5= \chi_1 \oplus 2 \chi_4 \oplus \chi_3 \oplus \chi'_3    \oplus 2 \chi_5 $$
 	where an index denotes the dimension of a character. By Corollary 6 we have  $$ \frac{5}{60} \; \equiv \;  c(A_5 ) \; \leq \; \frac{1}{25}( 1 + 2/4 + 2/3 + 2/5 ) \; = \; \frac{60 + 94}{ 25 \cdot 60 } $$
 \end{example}
\noindent 
 We will end this short note with a question.  Theorem 1 provides necessary conditions for a virtual character  (2.2)  to be non-negative. It is highly probable that these conditions are also sufficient. At this point, we can say however, that inequality (2.6) by itself is not enough to guarantee non-negativity of the relevant virtual character.
   For example, consider the sum of all irreducible characters
 \begin{align}
 	\mathfrak{\gamma} =	\psi_1 + \cdots +  \psi_k \nonumber
 \end{align}
 If we suppose that $ \gamma$ is non-negative (cf. Remark 1) then (2.5) and (1.2) yield

\begin{align}
\frac{ n_1 + \cdots +  n_k }{ n_1^2 + \cdots +  n_k^2 } \;  \leq \; 	
	\frac{ (1/n_1^2) \Re(\psi_1(g)) + \cdots + (1/n_q^{2})  \Re(\psi_k(g)) } 
	{k}
\end{align}
\noindent and (2.6) boils down to a well known inequality between harmonic and contraharmonic means (cf. e.g \cite{contraharmonic})  
\begin{align}
\frac{ n_1 + \cdots +  n_k }{ n_1^2 + \cdots +  n_k^2 }  \; \leq \;  \frac{(1/n_1) + \cdots + (1/n_k) }{ k}
\end{align}
\noindent On the other hand it was mentioned already (Remark 1)  that Gelfand character is not necessarily non-negative (\cite{Jorge})

\end{document}